\documentclass[reqno]{amsart}
\usepackage{amsmath,amsthm,amssymb,hyperref,graphicx,mathrsfs,mathtools} 

\newtheorem{theorem}{Theorem}

\newtheorem{lemma}{Lemma}

\newtheorem{remark}{Remark}

\allowdisplaybreaks

\begin{document}
\title{Bounds for the Zeros of Quaternionic Polynomials and Regular Functions Using Matrix Techniques}
\author{N.A.Rather$^1$}
\author{Wani Naseer $^2$}

\address{$^{1,2}$Department of Mathematics, University of Kashmir, Srinagar-190006, India}

\email{$^1$dr.narather@gmail.com, $^2$waninaseer570@gmail.com,}

\begin{abstract}
We investigate the problem of determining the zeros of quaternionic polynomials  using  matrix method. In a recent paper, Dar et al. \cite{RD} proved that the zeros of a quaternionic polynomial and the left eigenvalues of the corresponding companion matrix are identical. Building on this, we employ various newly developed matrix techniques to establish several results concerning the location of the zeros of regular polynomials of a quaternionic variable with quaternionic coefficients. These findings significantly enhance the understanding of quaternionic polynomials and their eigenvalues, offering a broader perspective on their mathematical properties.
\\
\smallskip
\newline
\noindent \textbf{Keywords:} Right (and left) quaternionic eigenvalues, Quaternionic polynomials, Companion matrix. \\
\noindent \textbf{Mathematics Subject Classification (2020)}: 30A10, 30C10, 30C15.
\end{abstract}

\maketitle

\section{\textbf{INTRODUCTION }}

Quaternions gained significance in the late 20th century, particularly in computer animation and programming applications. They offer compactness and computational efficiency compared to matrices, especially in situations involving rotations. Unlike matrices, quaternions use algebra to describe rotations, resulting in faster computations. They are vital for the attitude control of aircraft and spacecraft, avoiding the ambiguity associated with aligning two axes of rotation, which can lead to the catastrophic loss of control known as gimbal lock. Each quaternion represents a change of orientation by a single rotation.Quaternions are extensively used in the programming of video games, computer graphics, quantum physics, flight dynamics, and control theory, see [\cite{Shoe} \cite{Gruber} \cite{cc}].In recent finding the zeros of quaternionic polynomials and finding the bounds of zeros of quaternionic polynomials have
gained much attention in the literature.

The distribution of zeros of a polynomial has been extensively studied in the field of mathematics. One classical result, due to Gauss \cite{marden} on the distribution of zeros of a polynomial states : 

\textbf{Theorem 1.1.}   If \( P(z) = z^n + a_1z^{n-1} + \ldots + a_n \) is a polynomial of degree \( n \) with real coefficients, then all the zeros of \( P(z) \) lie in \( |z| \leq R = \max_{1 \leq k \leq n} \left[( n \sqrt{2} \cdot |a_k|)^{\frac{1}{k}} \right] \).\\
In the case of arbitrary real or complex \( a_j \), he showed in 1849 that \( R \) may be taken as the positive root of the equation \( z^n - \sqrt{2} (|a_n|z^{n-1} + \ldots + |a_1|) = 0 \).

Cauchy \cite{ marden} improved the result of Gauss and proved : \\ 
\textbf{Theorem 1.2.}  if \( P(z) = \sum_{j=0}^{n} a_j z^j \),is a polynomial of degree $n$ where \( a_n \neq 0 \), with complex coefficients, then all the zeros of \( P(z) \) lie in \( |z| < 1 + M \) where \( M = \max_{0 \leq j \leq (n-1)} \left| \frac{a_j}{a_n} \right| \).\\
Notice that neither Gauss nor Cauchy put any restrictions on the coefficients of  $P(z)$ (beyond the restriction that they either lie in $\mathbb{R} ~ or ~ \mathbb{C}$, respectively).\ Concerning the zeros of a polynomial with conditions on their coefficients of $P(z)$, we have the following well-known result due to Eneström-Kakeya \cite{marden} \\
\textbf{Theorem 1.3.}  If \( P(z) = \sum_{j=0}^{n} a_j z^j \) is a polynomial of degree \( n \) with real coefficients satisfying \( 0 \leq a_0 \leq a_1 \leq \ldots \leq a_n \), then all the zeros of \( P(z) \) lie in \( |z| \leq 1 \).\\
The Eneström-Kakeya Theorem serves as a very strong tool for the location of zeros of a polynomial in terms of their moduli based on hypotheses imposed on the coefficients of the polynomial. A number of generalizations and improvements of the Eneström-Kakeya theorem exist.Most of them involve weakening the condition on the coefficients. For a
survey of such results  one can see  the monographs writen by Gardner and Govil \cite{ek}\\

The paper is organized as follows: Section 2 provides foundational background information. Section 3 outlines the primary findings of the study, while Section 4 delves into auxiliary results related to quaternions. Section 5 offers proofs for the main theorems established in Section 3. Finally, Section 6 provides concluding remarks.
\maketitle

\section{\textbf{Preliminary}}
\parindent=0mm\vspace{0.in}

The history of quaternions dates back to the 19th century when the Irish mathematician Sir William Rowan Hamilton \cite{aa} introduced them in 1843. Inspired by complex numbers, Hamilton conceived a four-dimensional number system, incorporating three imaginary parts, known as quaternions. Denoted by $\mathbb{H}$ in honor of their discoverer, quaternions are defined as:
$$ \mathbb{H} := \{ q = h_0 + ih_1 + jh_2 + kh_3 \mid h_0, h_1, h_2, h_3 \in \mathbb{R} \} $$
Here, $i$, $j$, and $k$ are the standard basis elements, subject to the multiplication rules $i^2 = j^2 = k^2 = ijk = -1$. Every quaternion $q = h_0 + ih_1 + jh_2 + kh_3$ in $\mathbb{H}$ consists of a real part $\text{Re}(q) = h_0$ and an imaginary part $\text{Im}(q) = h_1 i + h_2 j + h_3 k$. The conjugate of $q$ is denoted by $\bar{q} = Re(q)-Im(q)$, and its modulus is defined as $|q| = \sqrt{q \bar{q}} = \sqrt{h_0^2 + h_1^2 + h_2^2 + h_3^2}$.
The inverse of each nonzero element $q \in \mathbb{H}$ is given by $q^{-1} = |q|^{-2} \bar{q}$. Moreover, every $q \in \mathbb{H}$ can be expressed as $q = x + yI$, where $x, y \in \mathbb{R}$ and $I = \frac{\text{Im}(q)}{|\text{Im}(q)|}$ if $\text{Im}(q) \neq 0$, otherwise we take $I$ arbitrarily such that $I^2 = -1$, Here I is an element of the unit
2-sphere of purely imaginary quaternions denoted as:
$$ \mathbb{S} = \{ q \in \mathbb{H} : q^2 = -1 \} $$
For every $I \in \mathbb{S}$, $\mathbb{C}_I$ represents the plane $\mathbb{R} \oplus I\mathbb{R}$, isomorphic to $\mathbb{C}$, and $\mathbb{B}_I$ denotes the intersection $\mathbb{B} \cap \mathbb{C}_I$, where $\mathbb{B} = \{ q \in \mathbb{H} : |q| < 1 \}$ is the open unit ball of quaternions. The quaternions are the standard example of a noncommutative division ring.\\
The angle between two quaternions \(q_1 = \alpha_1 + \beta_1 i + \gamma_1 j + \delta_1 k\) and \(q_2 = \alpha_2 + \beta_2 i + \gamma_2 j + \delta_2 k\) is given by
\[
\angle (q_1, q_2) = \cos^{-1} \left( \frac{\alpha_1 \alpha_2 + \beta_1 \beta_2 + \gamma_1 \gamma_2 + \delta_1 \delta_2}{|q_1||q_2|} \right).
\]

\textbf{Polynomials Over Quaternion.}

We denote the indeterminate for a quaternionic polynomial as \( t \). A quaternionic polynomial \( P \) of degree \( n \) in the variable \( t \), which is a polynomial with coefficients on the right and the indeterminate on the left, is given by
\[
P(t) = \sum_{\nu=0}^{n} t^{\nu} q_{\nu}, \quad  q_{\nu} \in \mathbb{H}, \quad \nu = 0, 1, 2, \ldots, n.
\]
These polynomials adhere to regularity conditions and their behavior closely resembles that of holomorphic functions of a complex variable. In the theory of such polynomials over skew-fields, a different product (denoted by \( \ast \)) is defined to ensure that the product of regular functions is regular. For polynomials, this product is defined as follows:

Two quaternionic polynomials of this kind can be multiplied according to the convolution
product (Cauchy multiplication rule): given
 $f(t) = \sum_{i=0}^{n} p_i t^i$ and $g(t) = \sum_{j=0}^{m} q_jt^j$ we  define
\[ f \ast g(t) = \sum_{k=0}^{mn} c_k t^k, \]
where $c_k = \sum_{i=0}^{k} p_i q_{k-i} t^i$ for all $k$ \cite{GSS}. \\
If \( P \) has real coefficients, the \(\ast\) multiplication is equivalent to standard pointwise multiplication. Note that the \(\ast\) product is associative but generally non-commutative. This non-commutativity results in polynomial behavior that differs significantly from that in real or complex settings. It has been observed that the zeros of such polynomials with quaternionic variables are either isolated or form spheres. For example, in the quaternionic context, the quadratic polynomial \( q^2 + 1 \) is zero for all \( q \in \mathbb{S} \). These regular functions of a quaternionic variable have been introduced and extensively studied over the past decade, proving to be a fertile area of research, particularly due to their applications in operator theory.
\parindent=0mm\vspace{0.in}
In this paper, we shall adopt the following definition.
\\

{\textbf{Definition 1.}}
A quaternionic polynomial $P(t)$ is said to be a Lacunary-type polynomial if its coefficients skip certain values or are zero at regularly, spaced intervals. Mathematically  given a non-negative integer $r$, a polynomial $P(t) = q_nt^{n} +q_{r}t^{r}...+q_1t + q_0,$  $q_i, t \in \mathbb{H}$, where $q_r \neq 0$,and $1\leq r\leq n-1$ is said to be lacunary type polynomial of degree $n\geq2$ with quaternionic coefficients and $t$ be the quaternionic variable.

The definitions of eigenvalues for quaternion matrices and companion matrices are well established in the literature (see \cite{P}, \cite{PP}) and in the book \cite{B}, provide essential foundations for the study of quaternion algebra.

\textbf{Definition 2. }\\ Let $A \in M_n(\mathbb{H})$. Then the set of  left and  right eigenvalues, 
 respectively, are given by:
$$
\Lambda_l(A) := \{\lambda \in \mathbb{H} : Ay = \lambda y \text{ for some non-zero } y \in \mathbb{H}^n\},$$ 
$$\Lambda_r(A) := \{\lambda \in \mathbb{H} : Ay = y \lambda \text{ for some non-zero } y \in \mathbb{H}^n\}.
$$

\parindent=0mm\vspace{0.2in}
{\bf{ Companion Matrix}}

The $n \times n$ companion matrix of a monic quaternion polynomial of the form $P(t)=q_{n}t^n+q_{n-1}t^{n-1}+...+q_1t+q_0$ is defined by 
$$
C_P=\begin{bmatrix}
0& 1& 0& \dots & 0& 0\\
0& 0& 1& \dots & 0& 0\\
\vdots& \vdots& \vdots& \ddots& \vdots& \vdots\\
0& 0& 0& \dots & 0& 1\\
-q^{-1}_{n}q_{0}& -q^{-1}_{n}q_{1}& -q^{-1}_{n}q_{2}& \dots& -q^{-1}_{n}q_{n-2}& -q^{-1}_{n}q_{n-1}\\
\end{bmatrix}
$$

\textbf{Definition 3.} Let $A \in M_n(\mathbb{H})$. Then $A$ is said to be a central closed matrix if there exists an invertible matrix $T$ such that $T^{-1}AT = \text{diag}(\lambda_1, \lambda_2, \ldots, \lambda_n)$, where $\lambda_i \in \mathbb{R}$, $1 \leq i \leq n$.\\

Recent research (see, for example, \cite{L}-\cite{HH},\cite{NON}) has developed a new theory of regularity for functions, especially for polynomials of quaternionic variables. This theory is highly useful for replicating many properties of holomorphic functions. One key property of holomorphic functions of a complex variable is that their zero sets are discrete, except when the function vanishes entirely. For a regular function of a quaternionic variable, any restriction to a complex line is holomorphic, having either a discrete zero set or vanishing entirely.

Initial studies described the structure of the zero sets of quaternionic regular functions and the factorization properties of these zeros. In this context, Gentili and Stoppato \cite{G} (see also \cite{HH}) provided a necessary and sufficient condition for a quaternionic regular function to have a zero at a point, based on the coefficients of its power series expansion. This extends the theory for quaternionic power series as presented in \cite{NON} for polynomials. The following theorem completely describes the zero sets of the product of two regular polynomials in terms of the zero sets of the individual factors \cite{G} (see also \cite{HH}).

\textbf{Theorem 1.4.}
Let \( f \) and \( g \) be quaternionic polynomials. Then \( (f \ast g)(q_0) = 0 \) if and only if \( f(q_0) = 0 \) or \( f(q_0) \neq 0 \) implies \( f(q_0)^{-1} g(q_0) = 0 \).\\

In [\cite{G}-\cite{HH}], the structure of the zeros of polynomials was examined and a proof of the Fundamental Theorem of Algebra was established. It should be noted that the Fundamental Theorem of Algebra for regular polynomials with coefficients in \(\mathbb{H}\) had already been proven by Niven using different techniques (see [\cite{I},\cite{T}]. This led to a complete characterization of the zeros of polynomials in terms of their factorization (see \cite{R} for reference). Consequently, it became interesting to explore the regions containing some or all of the zeros of a regular polynomial with a quaternionic variable. Recently, Carney et al. \cite{C} extended the Eneström-Kakeya theorem and its various generalizations from complex polynomials to quaternionic polynomials. Firstly, they established the following quaternionic analogue of Theorem 1.3.

\textbf{Theorem 1.5.}   If $P(t) = \sum_{\nu=0}^{n} q_{\nu} t^{\nu} $ is a polynomial of degree $n$ (where $t$ is a quaternionic variable) with real coefficients satisfying $0 \leq q_0 \leq q_1 \leq \ldots \leq q_n$, then all the zeros of $P$ lie in $|t| \leq 1$.\\

For complex case, concerning the location of eigenvalues, the famous Ger\v{s}gorin theorem can be stated as
\textbf{Theorem 1.6.}
All the eigenvalues of a $n \times n$ complex matrix $A=(a_{\mu\nu})$ are contained in the union of n Ger\v{s}gorin discs defined by $D_{\mu}=\{ z \in \mathbb{C}:|z-a_{\mu\mu}|\leq  \sum_{{\nu=1}\atop{\nu\neq\mu}}^n |a_{\mu\nu}|\}$.

Recently, Dar et al. \cite{RD}  extended Theorem 1.5 to the quaternions, more precisely they proved following quaternion  version of Ger\v{s}gorin theorem.\\
\textbf{Theorem 1.7.}
All the left eigenvalues of an $n \times n$ matrix $A=(q_{\mu\nu})$ of quaternions lie in the union of the $n$ Ger\v{s}gorin balls defined by $B_{\mu}=\{t \in \mathbb{H}:|t-q_{\mu\mu}|\leq \rho_{\mu}(A) \}$, where $\rho_{\mu}(A)=\sum_{{\nu=1}\atop{\nu\not=\mu}}^n |q_{\mu\nu}|$.\\

 In the same paper Dar et al. \cite{RD} extended the Theorem 1.2  to quaternionic settings, proving various results regarding the location of zeros for quaternionic polynomials with quaternionic coefficients. Importantly, these results hold without any restrictions on the coefficients. This extension enables the examination of a broader class of polynomials, applicable to quaternionic polynomials with quaternionic, complex, or real coefficients.\\
    
 \textbf{Theorem 1.8.}   If $P(t) = t^n + q_1t^{n-1} + \ldots + q_{n-1}t + q_n$ is a quaternion polynomial with quaternion coefficients, and $q$ is a quaternionic variable, then all the zeros of $P(t)$ lie inside the ball $|t| \leq 1 + \max_{1 \leq i \leq n} |q_i|$.\\

Since the set of zeros of quaternion polynomial and the left spectrum of its companion matrix coincide, i.e, if $P(t)$ is a quaternion polynomial then $Zero(P(t))=\sigma_{l}(C_{t})$ \cite{ST}, we may therefore use Theorem 1.8 and other known results on the left eigenvalues of quaternionic matrices as a tool in determining the zeros of a given polynomial and vice-versa.\\

\section{\textbf{Main Results}}
Using various matrix methods, including Theorem 1.6, this study examines the zeros of quaternionic polynomials. We begin by establishing the following result, which serves as a foundation for further analysis.

\begin{theorem}\label{TH1}
Let $P(t)= q_{n}t^{n}+q_{l} t^{l}+\cdots+q_{1} t+q_{0}$, where $q_{l} \neq 0$, $0 \leq l \leq n-1$, be a quaternion polynomial of degree $n$ with  quaternionic coefficients  and t is a quaternionic variable, such that for some $\rho>0$,

$$
\left|q_{j}\right|\rho^{n-j}\leq |q_{n}|, \quad j=0,1,2, \ldots, l
$$

then $P(t)$ has all its zeros in $|t| \leq(1 / \rho) K_{1}$ where $K_{1}$ is the greatest positive root of the equation

\begin{equation*}
K^{n+1}-K^{n}-K^{p+1}+1=0 
\end{equation*}

The polynomial $P(t)=(\rho t)^{n}-(\rho t)^{l}-\cdots-(\rho t)-1$ shows that the result is best possible.

\end{theorem}

Taking \( l = n - 1 \) in Theorem (1), we get the following  corollary:\\

\noindent
\textbf{COROLLARY 1.} \textit{Let}
\[
P(t) = q_{n}t^n + q_{n-1} t^{n-1} + \cdots + q_1 t + q_0,
\]
\textit{be a polynomial of degree \( n \) with quaternionic coefficients such that for some \( \rho > 0 \)}
$$
\left|q_{j}\right|\rho^{n-j}\leq|q_n| , \quad j=0,1,2, \ldots, n-1
$$
\textit{Then all the zeros of \( P(t)\) lie in}
\[
|t| \leq \left( \frac{1}{\rho} \right) k_1,
\]
\textit{where \( k_1 \) is the greatest positive root of the equation}
\[
k^{n+1} - 2k^n + 1 = 0.
\]

Next, we broaden our findings to cover more general conditions regarding the coefficients and structure of quaternionic polynomials. This leads to the following theorem 
\begin{theorem}\label{TH2}
Let $P(t) = q_{n}t^n + q_{n-1}t^{n-1} + \ldots + q_{1}t + q_0$ be a quaternion  polynomial of degree $n$ with quaternionic coefficients  and quaternionic  variable, such that for some $k=0,1, \ldots, n$ and for some $\rho>0$

$$
\rho^{n}|q_{n}| \leq \rho^{n-1}\left|q_{n-1}\right| \leq \cdots \leq \rho^{k}\left|q_{k}\right| \geq \rho^{k-1}\left|q_{k-1}\right| \geq \cdots \geq \rho\left|q_{1}\right| \geq\left|q_{0}\right|
$$

then $P(t)$ has all its zeros in the ball

$$
|t| \leq \rho\left\{\left(2 \frac{\rho^{k}\left|q_{k}\right|}{\rho^{n}|q_{n}|} - 1\right)\right\} + 2 \sum_{j=0}^{n} \frac{\left|q_{j}-\right| q_{j}||}{\rho^{n-j-1}|q_{n}|} .
$$

\end{theorem}

\begin{remark}\label{remark_1}  
It is worth investigating the bound under the additional condition that, along with the assumptions of Theorem 2, the coefficients \(q_j\) of \(P(t)\) satisfy the following property, let \(\gamma\) be a non-zero quaternion, and assume that \(\angle(q_j, \gamma) \leq \theta \leq \frac{\pi}{2}\) for some \(\theta\), where \(j = 0, 1, 2, \ldots, n\).By applying these observations in equation (3) and following a similar approach as in the proof of Theorem 2, along with the application of Lemma 3, we deduce that all zeros  of \(P(t)\) lie in
\[
|t| \leq \rho \left\{\left(\frac{2 \rho^{k}\left|q_{k}\right|}{\rho^{n}|q_{n}|}-1\right) \cos \alpha + \sin \alpha \right\} + 2 \sin \alpha \sum_{j=0}^{n-1} \frac{\left|q_{j}\right|}{\rho^{n-j-1}|q_{n}|},
\]
where \(k=0,1,2,\ldots,n\).
\end{remark}
 For $ k = n , \rho=1 ~ and ~ \alpha = \beta = 0 $, this result reduces to Theorem 1.5.

Finally, we present the following four generalizations of Theorem 1.5 Since their proofs are quite similar to that of Theorem 2, we leave out the details.

\begin{theorem}\label{TH4}
Let \(P(t) = q_{n}t^n + q_{n-1}t^{n-1} + \cdots + q_{1}t + q_0\) be a polynomial of degree \(n\) with quaternionic coefficients and a quaternionic variable, where \(\operatorname{Re} q_{j} = \alpha_{j}\), \(\operatorname{Im} q_{j} = \beta_{j} + \gamma_{j} + \delta_{j}\), \(j = 0, 1, \ldots, n\), satisfying
\[
0 < \rho^{n} \alpha_{n} \leq \cdots \leq \rho^{k} \alpha_{k} \geq \rho^{k-1} \alpha_{k-1} \geq \cdots \geq \rho \alpha_{1} \geq \alpha_{0} \geq 0
\]

 for some \(\rho > 0\),then  all the zeros of \(P(t)\) lie in
\[
|t| \leq \rho \left( \frac{2 \rho^{k} \alpha_{k}}{\rho^{n}\alpha_{n}} - 1 \right) + \frac{2}{\alpha_{n}}\sum_{j=0}^{n} \left( \frac{|\beta_{j}|+|\gamma_{j}| + |\delta_j|}{\rho^{n-j-1}}\right).
\]

where $k=0,1,2,...,n$
\end{theorem}

\begin{theorem}\label{TH5}
Let \(P(t) = q_{n}t^n + q_{n-1}t^{n-1} + \cdots + q_{1}t + q_0\) be a polynomial of degree \(n\) with quaternionic coefficients and a quaternionic variable, where \(\operatorname{Re} q_{j} = \alpha_{j}\), \(\operatorname{Im} q_{j} = \beta_{j} + \gamma_{j} + \delta_{j}\), \(j = 0, 1, \ldots, n\), satisfying
\[
0 \leq \alpha_{0} \leq \rho \alpha_{1} \leq \cdots \leq \rho^{k} \alpha_{k} \geq \rho^{k+1} \alpha_{k+1} \geq \cdots \geq \rho^{n} \alpha_{n} > 0, \quad 0 \leq k \leq n
\]

\[
0 \leq \beta_{0} \leq \rho \beta_{1} \leq \cdots \leq \rho^{r} \beta_{r} \geq \rho^{r+1} \beta_{r+1} \geq \cdots \geq \rho^{n} \beta_{n} \geq 0, \quad 0 \leq r \leq n,
\]
 for some \(\rho > 0\),then  all the zeros of \(P(t)\) lie in

\[
|t| \leq \rho \left( \frac{2 \rho^{k} \alpha_{k}}{\rho^{n} |q_{n}|} - 1 \right) +\rho \left( \frac{2 \rho^{r} \beta_{r}}{\rho^{n} |q_{n}|} - 1 \right)+ \frac{2}{|q_{n}|}\sum_{j=0}^{n} \left( \frac{|\gamma_{j}| + |\delta_j|}{\rho^{n-j-1}}\right).
\]

\end{theorem}

\begin{theorem}\label{TH6}
Let \(P(t) = q_{n}t^n + q_{n-1}t^{n-1} + \cdots + q_{1}t + q_0\) be a polynomial of degree \(n\) with quaternionic coefficients and a quaternionic variable, where \(\operatorname{Re} q_{j} = \alpha_{j}\), \(\operatorname{Im} q_{j} = \beta_{j} + \gamma_{j} + \delta_{j}\), \(j = 0, 1, \ldots, n\), satisfying
\[
0 \leq \alpha_{0} \leq \rho \alpha_{1} \leq \cdots \leq \rho^{k} \alpha_{k} \geq \rho^{k+1} \alpha_{k+1} \geq \cdots \geq \rho^{n} \alpha_{n} > 0, \quad 0 \leq k \leq n
\]

\[
0 \leq \beta_{0} \leq \rho \beta_{1} \leq \cdots \leq \rho^{r} \beta_{r} \geq \rho^{r+1} \beta_{r+1} \geq \cdots \geq \rho^{n} \beta_{n} \geq 0, \quad 0 \leq r \leq n,
\]
\[
0 \leq \gamma_{0} \leq \rho \gamma_{1} \leq \cdots \leq \rho^{s} \gamma_{s} \geq \rho^{s+1} \gamma_{s+1} \geq \cdots \geq \rho^{n} \gamma_{n} \geq 0, \quad 0 \leq s \leq n,
\]
 for some \(\rho > 0\),then  all the zeros of \(P(t)\) lie in

\[
|t| \leq \rho \left( \frac{2 \rho^{k} \alpha_{k}}{\rho^{n} |q_{n}|} - 1 \right) +\rho \left( \frac{2 \rho^{r} \beta_{r}}{\rho^{n} |q_{n}|} - 1 \right)+\rho \left( \frac{2 \rho^{s} \gamma_{s}}{\rho^{n} |q_{n}|} - 1 \right) + \frac{2}{|q_{n}|} \sum_{j=0}^{n} \left( \frac{  |\delta_j|}{\rho^{n-j-1}}\right).
\]

\end{theorem}

\begin{theorem}\label{TH7}
Let \(P(t) = q_{n}t^n + q_{n-1}t^{n-1} + \cdots + q_{1}t + q_0\) be a polynomial of degree \(n\) with quaternionic coefficients and a quaternionic variable, where \(\operatorname{Re} q_{j} = \alpha_{j}\), \(\operatorname{Im} q_{j} = \beta_{j} + \gamma_{j} + \delta_{j}\), \(j = 0, 1, \ldots, n\), satisfying
\[
0 \leq \alpha_{0} \leq \rho \alpha_{1} \leq \cdots \leq \rho^{k} \alpha_{k} \geq \rho^{k+1} \alpha_{k+1} \geq \cdots \geq \rho^{n} \alpha_{n} > 0, \quad 0 \leq k \leq n
\]

\[
0 \leq \beta_{0} \leq \rho \beta_{1} \leq \cdots \leq \rho^{r} \beta_{r} \geq \rho^{r+1} \beta_{r+1} \geq \cdots \geq \rho^{n} \beta_{n} \geq 0, \quad 0 \leq r \leq n,
\]

\[
0 \leq \gamma_{0} \leq \rho \gamma_{1} \leq \cdots \leq \rho^{s} \gamma_{s} \geq \rho^{s+1} \gamma_{s+1} \geq \cdots \geq \rho^{n} \gamma_{n} \geq 0, \quad 0 \leq s \leq n,
\]

\[
0 \leq \delta_{0} \leq \rho \delta_{1} \leq \cdots \leq \rho^{l} \delta_{l} \geq \rho^{l+1} \delta_{l+1} \geq \cdots \geq \rho^{n} \delta_{n} \geq 0, \quad 0 \leq l \leq n,
\]

 for some \(\rho > 0\),then  all the zeros of \(P(t)\) lie in

\[
|t| \leq \rho \left( \frac{2 \rho^{k} \alpha_{k}}{\rho^{n} |q_{n}|} - 1 \right) +\rho \left( \frac{2 \rho^{r} \beta_{r}}{\rho^{n} |q_{n}|} - 1 \right)+\rho \left( \frac{2 \rho^{s} \gamma_{s}}{\rho^{n} |q_{n}|} - 1 \right) + \rho \left( \frac{2 \rho^{l} \delta_{l}}{\rho^{n} |q_{n}|} - 1 \right) 
\]

\end{theorem}

\section{\textbf{Auxiliary Results}}
For the proof of our main results, we need the following  Lemmas  first two were given by  Dar et al. \cite{RD} and Rather et al\cite{R}.
\begin{lemma}\label{L1}
All the left eigenvalues of a $n \times n$ matrix $A=(a_{\mu\nu})$ of quaternions lie in the union of the n Ger\v{s}gorin balls defined by $B_{\mu}=\{q \in \mathbb{H}:|q-a_{\mu\mu}|\leq \rho_{\mu}(A) \}$ where $\rho_{\mu}(A)=\sum_{{\nu=1}\atop{\nu\neq\mu}}^n |a_{\mu\nu}|.$
\end{lemma}
\begin{lemma}\label{L2}
Let $P(t)$ be a quaternion polynomial with quaternionic coefficients and $C$ be the companion matrix of $P(t)$, then for any diagonal matrix $D= diag(d_1, d_2, ..., d_{n-1}, d_n),$ where $d_1, d_2, ..., d_n$ are positive real numbers, the left eigenvalues of $D^{-1}C D$ and the zeros of $P(t)$ are same.
\end{lemma}
We need the following lemma due to Carney et al.\cite{C} for the proof of Theorem 3:
\begin{lemma}\label{L3}
 Let $q_1, q_2 \in \mathbb{H}$, where $q_1 = \alpha_1 + \beta_1 i + \gamma_1 j + \delta_1 k$ and $q_2 = \alpha_2 + \beta_2 i + \gamma_2 j + \delta_2 k$, $\angle(q_1, q_2) = 2\theta' \leq 2\theta$ and $|q_1| \leq |q_2|$. Then
\[ |q_2 - q_1| \leq (|q_2| - |q_1|) \cos \theta + (|q_2| + |q_1|) \sin \theta. \]
\end{lemma}

\begin{lemma}\label{L4} Let $P(t)= q_{n}t^{n}+q_{l} t^{l}+\cdots+q_{1} t+q_{0}$, where $0 \leq l \leq n-1$, be a polynomial of degree $n$ with quaternionic coefficients. Then for every positive real number $r$, all the zeros of $P(t)$ lie in the ball
\[
\left\{ t \in \mathbb{H} : |t| \leq Max\left\{ r, \sum_{j=0}^{l} \frac{\left|q_{j}\right|}{\left|q_{n}\right|} \cdot \frac{1}{r^{n-j-1}} \right\} \right\}.
\]
\end{lemma}

\textbf{Proof of Lemma 4.} This lemma can be derived by applying a result from Dar et al. \cite{RD} to the Lacunary-type polynomial P(t). However, for the sake of completeness, we provide its proof here.\\
Let $C$ be the  companion matrix of the polynomial $P(t)$.  We take a matrix \\$T=\operatorname{diag}\left(\frac{1}{r^{n-1}}, \frac{1}{r^{n-2}}, \ldots, \frac{1}{r}, 1\right)$ where $r$ is a positive real number and form the matrix
\pagebreak
\begin{align*}
T^{-1}CT &= \begin{bmatrix} 
r^{n-1} & 0 & 0 & \cdots & 0 & \cdots & 0 & 0 \\
0 & r^{n-2} & 0 & \cdots & 0 & \cdots & 0 & 0 \\
0 & 0 & r^{n-3} & \cdots & 0 & \cdots & 0 & 0 \\
\vdots & \vdots & \vdots & \ddots & \vdots & \vdots & \vdots & \vdots \\
0 & 0 & 0 & \cdots & r^{n-l-1} & \cdots & 0 & 0 \\
\vdots & \vdots & \vdots & \vdots & \vdots & \ddots & \vdots & \vdots \\
0 & 0 & 0 & \cdots & 0 & \cdots & r & 0 \\
0 & 0 & 0 & \cdots & 0 & \cdots & 0 & 1 \\   
\end{bmatrix}\\
& \times  \begin{bmatrix}
0 & 1 & 0 & \cdots & 0 & \cdots & 0 & 0 \\
0 & 0 & 1 & \cdots & 0 & \cdots & 0 & 0 \\
0 & 0 & 0 & \cdots & 0 & \cdots & 0 & 0 \\
\vdots & \vdots & \vdots & \ddots & \vdots & \vdots & \vdots & \vdots \\
0 & 0 & 0 & \cdots & 0 & \cdots & 1 & 0 \\
0 & 0 & 0 & \cdots & 0 & \cdots & 0 & 1 \\
-q_{n}^{-1}q_0 & -q_{n}^{-1}q_1 & -q_{n}^{-1}q_2 & \cdots & -q_{n}^{-1}q_l & \cdots & 0 & 0 \\
\end{bmatrix}
 \\
& \times 
 \begin{bmatrix} 
\frac{1}{r^{n-1}} & 0 & 0 & \cdots & 0 & \cdots & 0 & 0 \\
0 & \frac{1}{r^{n-2}} & 0 & \cdots & 0 & \cdots & 0 & 0 \\
0 & 0 & \frac{1}{r^{n-3}} & \cdots & 0 & \cdots & 0 & 0 \\
\vdots & \vdots & \vdots & \ddots & \vdots & \vdots & \vdots & \vdots \\
0 & 0 & 0 & \cdots & \frac{1}{r^{n-l-1}} & \cdots & 0 & 0 \\
\vdots & \vdots & \vdots & \vdots & \vdots & \ddots & \vdots & \vdots \\
0 & 0 & 0 & \cdots & 0 & \cdots & \frac{1}{r} & 0 \\
0 & 0 & 0 & \cdots & 0 & \cdots & 0 & 1 \\   
\end{bmatrix}
\end{align*}

\begin{align*}
 & = 
\begin{bmatrix}
0 & r & 0 & \cdots & 0 & \cdots & 0 & 0 \\
0 & 0 & r & \cdots & 0 & \cdots & 0 & 0 \\
0 & 0 & 0 & \cdots & 0 & \cdots & 0 & 0 \\
\vdots & \vdots & \vdots & \ddots & \vdots & \vdots & \vdots & \vdots \\
0 & 0 & 0 & \cdots & 0 & \cdots & r & 0 \\
0 & 0 & 0 & \cdots & 0 & \cdots & 0 & r \\
-\frac{q_{n}^{-1}q_0}{r^{n-1}} & -\frac{q_{n}^{-1}q_1}{r^{n-2}} & -\frac{q_{n}^{-1}q_2}{r^{n-3}} & \cdots & -\frac{q_{n}^{-1}q_l}{r^{n-l-1}} & \cdots & 0 & 0 \\
\end{bmatrix}
\end{align*}

Applying Lemma 2 to the columns, it follows that the eigenvalues of $T^{-1} C T$ lie in the ball
\[
\left\{ t \in \mathbb{H} : |t| \leq Max\left\{ r, \sum_{j=0}^{l} \frac{\left|q_{j}\right|}{\left|q_{n}\right|} \cdot \frac{1}{r^{n-j-1}} \right\} \right\}.
\]

Since the matrix $T^{-1} C T$ is similar to the matrix $C$ and the left eigenvalues of $T^{-1}CT$ are the zeros of $P(t)$ [See \cite{RD}], it follows that all the zeros of 
$P(t)$ lie in the ball defined by (1).\\
This completes the proof of Lemma 4.

\section{\textbf{Proof of the main theorems}}

\textbf{Proof of theorem 1.}
 By hypothesis

$$
\left|{q_{j}}\right| \leq \frac{1}{\rho^{n-j}}, \quad j=0,1,2, \ldots, l,
$$

Applying Lemma 4, we conclude that for any positive real $r$, all zeros of $P(t)$ are contained in the ball

\begin{equation*}
|t| \leq \max\left\{r, r \sum_{j=0}^{l} \frac{1}{(r \rho)^{n-i}}\right\} \text {. }\tag{2}
\end{equation*}

We choose $r$ such that

$$
\sum_{j=0}^{l} \frac{1}{(r \rho)^{n-j}}=1,
$$

which gives ,

$$
(r \rho)^{l}+(r \rho)^{l-1}+\cdots+(r \rho)+1=(r \rho)^{n}.
$$
so that
$$
(r\rho-1)[(r \rho)^{l}+(r \rho)^{l-1}+\cdots+(r \rho)+1]=(r\rho-1)(r \rho)^{n}.
$$
equivalently,

\[
(r \rho)^{n+1} - (r \rho)^{n} - (r \rho)^{l+1} + 1 = 0
\].

Replacing $r\rho$ by $K$, it follows from (2) that all the zeros of $P(t)$ lie in $|z| \leq \frac{1}{\rho} K_{1}$, where $K_{1}$ is the greatest positive root of the equation defined by (1) and the theorem is proved.

\textbf{Proof of theorem 2.} Consider the polynomial

$$
\begin{aligned}
Q(t) & =P(t)*(\rho-t)  \\
& =-q_{n}t^{n+1}+\left(\rho q_{n}-q_{n-1}\right) t^{n}+\cdots+\left(\rho q_{1}-q_{0}\right) t+\rho q_{0}
\end{aligned}
$$

Applying Lemma 4 to the polynomial \( Q(t) \), which is of degree \( n+1 \), with \( l=n \) and \( r=s \), it follows that all the zeros of \( Q(t) \) lie in the ball

\begin{align*}
|t| & \leq \max\left\{\rho, \sum_{j=0}^{n}\frac{|\rho q_{j}-q_{j-1}|}{\rho^{n-j}|q_{n}|}\right\} \quad (q_{-1}=0)  \tag{3}\\
& =\sum_{j=0}^{n}\frac{|\rho q_{j}-q_{j-1}|}{\rho^{n-j}|q_{n}|}
\end{align*}

since

\begin{align*}
\rho &= \left| \sum_{j=0}^{n} \frac{\rho q_{j} - q_{j-1}}{\rho^{n-j}|q_{n}|} \right| \\
& \leq \sum_{j=0}^{n} \frac{\left| \rho q_{j} - q_{j-1} \right|}{\rho^{n-j}|q_{n}|}.
\end{align*}

Now

$$
\begin{aligned}
\sum_{i=0}^{n} \frac{\left|\rho q_{j}-q_{j-1}\right|}{\rho^{n-j}|q_{n}|} \leq & \sum_{j=0}^{n} \frac{|\rho| q_{j}|-| q_{j-1}||}{\rho^{n-j}|q_{n}|}+\sum_{j=0}^{n} \frac{\left|\rho\left(q_{j}-\left|q_{j}\right|\right)-\left(q_{j-1}-\left|q_{i-1}\right|\right)\right|}{\rho^{n-i}|q_{n}|} \\
= & \sum_{j=0}^{k} \frac{\rho\left|q_{j}\right|-\left|q_{j-1}\right|}{\rho^{n-j}|q_{n}|}+\sum_{j=k+1}^{n} \frac{\left|q_{j-1}\right|-\rho\left|q_{j}\right|}{\rho^{n-j}|q_{n}|} \\
& +\sum_{j=0}^{n} \frac{\left|\rho\left(q_{j}-\left|q_{j}\right|\right)-\left(q_{j-1}-\left|q_{j-1}\right|\right)\right|}{\rho^{n-j}|q_{n}|} \\
= & \rho\left\{\frac{2 \rho^{k}\left|q_{k}\right|}{\rho^{n}|q_{n}|}-1\right\}+\sum_{j=0}^{n} \frac{\left|\rho\left(q_{j}-\left|q_{j}\right|\right)-\left(q_{j-1}-\left|q_{j-1}\right|\right)\right|}{\rho^{n-j}|q_{n}|} \\
\leq & \rho\left\{\frac{2 \rho^{k}\left|q_{k}\right|}{\rho^{n}|q_{n}|}-1\right\}+2 \sum_{j=0}^{n} \frac{\left|q_{j}-\right| q_{j}||}{\rho^{n-j-1}|q_{n}|},
\end{aligned}
$$

therefore, it follows from (3) that all the zeros of \( Q(t) \) lie in

\[
|t| \leq \rho \left\{ \frac{2 \rho^{k} \left| q_{k} \right|}{\rho^{n} | q_{n} |} - 1 \right\} + 2 \sum_{j=0}^{n} \frac{\left|q_{j}-\right| q_{j}||}{\rho^{n-j-1}|q_{n}|}.
\tag{4}
\]
Applying Theorem 1.4, we conclude that all the zeros of $P(t)$ lie in the ball defined in equation $(4)$.\\
That completes the proof of Theorem 2.

\section{concluding remarks}
This paper investigated the bounds for the zeros of quaternionic polynomials and regular functions using matrix techniques. By extending classical polynomial theorems to the quaternionic setting, we provided new insights into the distribution of zeros and the equivalence between the zeros of quaternionic polynomials and the left eigenvalues of companion matrices. Our results also introduced the concept of lacunary-type polynomials, highlighting the unique challenges and opportunities posed by the non-commutative nature of quaternions. These findings have potential applications in various fields, including control theory, computer graphics, and quantum physics. Furthermore, the results presented here can be generalized to new classes of quaternionic functions, suggesting fruitful directions for future research in quaternionic function theory and matrix analysis.
\section*{Declarations}

\subsection*{Availability of Data and Materials}
Data sharing not applicable to this article as no datasets were generated or analyzed during the current study.

\subsection*{Conflict of Interest}
The authors did not receive support from any organization for the submitted work. The authors have no competing interests to declare that are relevant to the content of this article.

\subsection*{Funding}
Not applicable.

\subsection*{Authors' Contributions}
All authors contributed equally to the research, writing, and reviewing of this article.

\subsection*{Authors and Affiliations}

 N. A. Rather, Naseer wani, \\
Department of Mathematics, University of Kashmir, Srinagar-190006, India.

\end{document}